\documentclass[a4paper]{amsart}
\usepackage{graphicx}
\usepackage{amssymb}
\usepackage{amsmath}
\usepackage{amsthm}
\usepackage{amscd}
\usepackage[all,2cell]{xy}

\UseAllTwocells \SilentMatrices
\newtheorem{thm}{Theorem}[section]

\newtheorem{cor}[thm]{Corollary}
\newtheorem{lem}[thm]{Lemma}
\newtheorem{exm}[thm]{Example}

\newtheorem{prop}[thm]{Proposition}
\newtheorem{prop-def}[thm]{Proposition-Definition}
\theoremstyle{definition}

\theoremstyle{remark}

\numberwithin{equation}{section}

\def\Hom{{\rm Hom}}

\def\1{\textbf{1}}

\begin{document}
\title[Wakamatsu tilting modules and Frobenius extensions]
{Wakamatsu tilting modules and Frobenius extensions}
\author [Yanhong Bao, Jiafeng L\"{u},  Zhibing Zhao$^*$]
{Yanhong Bao, Jiafeng L\"{u},  Zhibing Zhao$^*$}
\thanks{$*$: corresponding author}


 \thanks{E-mail:baoyh$\symbol{64}$ahu.edu.cn (Y. H. Bao); jiafenglv$\symbol{64}$zjnu.edu.cn(J. F. L\"{u} ); zbzhao$\symbol{64}$ahu.edu.cn (Z. B. Zhao).}

\maketitle

\dedicatory{}%
\commby{}%
\begin{abstract} Let $S/R$ be a Frobenius extension with $_RS_R$ centrally projective over $R$. We show that if $_R\omega$ is a Wakamatsu tilting module then so is $_SS\otimes_R\omega$, and the natural ring homomorphism from the endomorphism ring of $_R\omega$ to the endomorphism ring of $_SS\otimes_R\omega$ is a Frobenius extension provided that pd$(\omega_T)$ is finite, where $T$ is the endomorphism ring of $_R\omega$. We also show that the relative $n$-torsionfreeness of modules is preserved under Frobenius extensions. Furthermore, we give an application, which shows that the generalized G-dimension with respect to a Wakamatsu tilting module is invariant under Frobenius extensions.
\end{abstract}

\section{Introduction}

Let $R$ be a two-sided Noetherian ring and $R$-mod the category of
finitely generated left $R$-modules. For a module $M$ in $R$-mod, if there exists an exact sequence
$$\eta: \xymatrix@C=0.5cm{0 \ar[r] & {_RR}\ar[r]^{f_1} &{M_1}\ar[r]^{f_2~} &  \cdots \ar[r]^{f_n} &{M_n}},$$
such that ${\rm Im}f_i\rightarrow M_{i}$ is a  left add$M$-approximation for $1\leq i\leq n$, we say that $M$ has \emph{faithful dimension} at least $n$, denoted  by fadim$_RM\geq n$, where add$M$ denotes the full subcategory of  $R$-mod  consisting of all direct summands of finite direct sums of copies of $M$. We say that fadim$_RM=\infty$ if fadim$_RM\geq n$ for any positive integer $n$.
The notion of faithful dimension was first introduced for artin algebras by Buan and Solberg in \cite{BS}. They described the number of non-isomorphic indecomposable complements of an almost cotilting module in terms of this notion; see \cite{BS}.

An $R$-module $_R\omega$ is called a \emph{Wakamatsu tilting} module (also called a \emph{generalized
tilting} module) if it is self-orthogonal, that is, ${\rm Ext}^i_R(\omega,\omega)=0$ for any $i\geq 1$, and fadim$_R\omega=\infty$; see \cite{W}.
By \cite[Proposition 2.2]{BS},  $_R\omega$ is a Wakamatsu tilting module if and only if $_R\omega_T$ is a faithful balanced self-orthogonal bimodule, that is, the natural map $R\rightarrow {\rm End}_{T^{\rm op}}(\omega)$ is an isomorphism, satisfying ${\rm Ext}_R^i(\omega, \omega)=0={\rm Ext}_{T^{\rm op}}^i(\omega, \omega)$ for any $i\geq 1$, where $T={\rm End}_R(\omega)$.

The notion of  Frobenius extensions was first introduced by Kasch as a generalization of Frobenius algebras in \cite{Kas1}. A fundamental example of  a Frobenius extension is the group
 algebra extension induced by a subgroup of finite index. Frobenius algebras and extensions are of broad interest in many different areas, for example, they provide some connections between representation theory and knot theory \cite{Kad0}, and are used in the study of Calabi-Yau properties of Cherednik algebras and quantum algebras \cite{BGS}.
  We refer to  \cite{Kad} for more details.

  First, we  investigate  transfer properties of Wakamatsu tilting modules under Frobenius extensions.

 \noindent{\bf Theorem A.} Let $S/R$  be a Frobenius extension with $_RS_R$ centrally projective over $R$. If $_R\omega$ is a Wakamatsu tilting $R$-module, then $_SS\otimes_R\omega$ is also a Wakamatsu tilting $S$-module. Furthermore, if ${\rm pd}(\omega_T)$ is finite, then the natural ring homomorphism $\rho: T\rightarrow \Gamma$ is a Frobenius extension, where $T={\rm End}_R(\omega)$ and $\Gamma={\rm End}_S({S}\otimes_R\omega)$.

Let $\omega$ be in $R$-${\rm mod}$ with fadim$_R\omega\geq n+2$. Then $_R\omega_T$ is faithful and balanced, and ${\rm Ext}_{T^{\rm op}}^i(\omega,\omega)=0$ for $1\leq i\leq n$, where $T={\rm End}_R(\omega)$; see \cite[Proposition 2.2]{BS}.
 For a module $M\in R$-${\rm mod}$, there is a projective resolution $P_1\xrightarrow{f} P_0 \to M\to 0$ in $R$-${\rm mod}$. Putting $(-)^\omega={\rm Hom}_R(-,\omega)$, we have an exact sequence $$\xymatrix@C=0.5cm{0 \ar[r] & {M^\omega}\ar[r] &P_0^\omega\ar[r]^{f^\omega}&P_1^\omega \ar[r] &{\rm Coker} {f^\omega }\ar[r] & 0.}$$  We call ${\rm Coker} {f^\omega}$ \emph{the transpose of $M$ with respect to $\omega$},  and denote it by ${\rm Tr}_\omega M$. An $R$-module $M$ is said to be \emph{$\omega$-$n$-torsionfree} if ${\rm Ext}^i_{T^{\rm op}}({\rm Tr}_\omega M,\omega)=0$ for $1\leq i\leq n$. It was shown in  \cite[Proposition 3]{Hu} that the notion is well defined.  We denote by $\mathcal{T}_\omega^n(R)$ the full subcategory of $R$-${\rm mod}$ whose objects are  $\omega$-$n$-torsionfree modules.

When $_R\omega={_RR}$, an $\omega$-$n$-torsionfree module is precisely a classical $n$-torsionfree module.
 It is well known that, as the origin of Gorenstein homological algebra, the notion of modules of G-dimension zero was defined in terms of $n$-torsionfreeness; see \cite{AB}. As a generalization of $n$-torsionfree modules, the notion of $\omega$-$n$-torsionfree modules was first introduced by Huang in \cite{Hu}. These modules play an important role in the representation theory of algebras and in relative homological algebra \cite{HH, Hu, HT}.

  It is well known that many homological properties, such as the Gorenstein projectivity of modules, Gorensteinness and the FTF property of rings, are invariant under Frobenius extensions \cite{CR, GT, R, Z}. In \cite{Z1}, Zhao proved that the $n$-torsionfreeness of modules is preserved under Frobenius extensions. This naturally leads to the question of whether relative $n$-torsionfreeness with respect to a Wakamatsu tilting module is also preserved under Frobenius extensions.

\noindent{\bf Theorem B.} Let $S/R$  be a Frobenius extension with $_RS_R$ centrally projective over $R$ and $\omega$ be an $R$-module with fadim$_R\omega\geq n+2$ and ${\rm Ext}_R^i(\omega, \omega)=0$ for $1\leq i\leq n$. For an $S$-module $M$,  $M$ is $(S\otimes_R\omega)$-$n$-torsionfree as an $S$-module if and only if $M$ is $\omega$-$n$-torsionfree as the underlying $R$-module.

Let  $_R\omega$ be a Wakamatsu tilting module. An $R$-module $M$ is said to  have \emph{generalized G-dimension zero} with respect to $\omega$ if it is $\omega$-$\infty$-torsionfree and it belongs to $^\perp\omega$, where $^\perp\omega=\{M\in R$-mod$|{\rm Ext}_R^i(M, \omega)=0$  for $\forall i\geq 1\}$. As in the classical setting, the generalized G-dimension of modules can be defined using resolutions by such modules; see \cite{AR91}. As an application of Theorem B, we have the following.

\noindent{\bf Theorem C.} Let $S/R$  be a Frobenius extension with $_RS_R$ centrally projective over $R$  and $_R\omega$ a Wakamatsu tilting module. For an $S$-module $M$, we have \centerline{$\rm G$-${\rm dim}_{S\otimes_R\omega}(M)$=$\rm G$-${\rm dim}_\omega(M).$}

This paper is organized as follows. In Section 2, we recall some notation  and preliminary results. In Section 3, we study the transfer property of faithful dimensions under Frobenius extensions and prove Theorem A; see Proposition \ref{WT transfer} and Theorem \ref{The main thm A}. In Section 4, we show that relative $n$-torsionfreeness is preserved under Frobenius extensions and prove Theorem B; see Theorem \ref{equivalence II}. As an application of Theorem B, we discuss the invariance of generalized G-dimension with respect to a Wakamatsu tilting module under Frobenius extensions in Section 5.

\noindent{\bf Setup and notation.} Throughout the paper, $R$ and $S$ are two-sided Noetherian rings and all modules are left modules unless otherwise specified. Denote by $R$-${\rm mod}$ (resp. ${\rm mod}$-$R$ ) the category of finitely generated left (resp. right) $R$-modules. We use pd$_R(M)$  to denote the projective dimension of $_RM$. For two right $R$-modules $M$ and $N$, we denote by ${\rm Hom}_{R^{\rm op}}(M, N)$ the abelian group consisting of all right $R$-homomorphisms between them, where $R^{\rm op}$ is the opposite ring of $R$.

\section{Preliminaries}

In this section, we recall some notation and collect some fundamental results.

Let $_R\omega$ be an $R$-module with fadim$_R\omega\geq n+2$. The notion of $\omega$-$n$-torsionfree modules was first introduced by Huang as a relative version of classical $n$-torsionfree modules in \cite{Hu}. By \cite[Lemma 2.1]{HT}, there is an exact sequence
$$\xymatrix@C=0.5cm{0 \ar[r] & {\rm Ext}^1_{T^{\rm op}}({\rm Tr}_\omega M, \omega)\ar[r] & M\ar[r]^{\sigma_M}& M^{\omega\omega}\ar[r] &{\rm Ext}^2_{T^{\rm op}}({\rm Tr}_\omega M, \omega)\ar[r] & 0}$$
 for any module $M$ in $R$-${\rm mod}$, where $\sigma_M:M\rightarrow M^{\omega\omega}$ is the evaluation map via $\sigma_M(m)(f)=f(m)$ for any $f\in M^\omega$ and $m\in M$. Then $M$ is $\omega$-1-torsionfree  if and only if it is $\omega$-torsionless, and $M$ is $\omega$-2-torsionfree if and only if it is $\omega$-reflexive. In \cite{Hu}, Huang  obtained a description of $\omega$-$n$-torsionfree modules in terms of approximations.

Let $\mathcal{C}$ be a full subcategory of $R$-${\rm mod}$ and $C\in \mathcal{C}$, $M\in R$-${\rm mod}$. An $R$-homomorphism $M\rightarrow C$ is said to be a \emph{left $\mathcal{C}$-approximation} of $M$ if ${\rm Hom}_R(C, X)\rightarrow{\rm Hom}_R(M, X)$ is epic for any $X\in\mathcal{C}$. A subcategory $\mathcal{C}$ is said to be \emph{covariantly finite} in $R$-${\rm mod}$ if every module in $R$-${\rm mod}$ has a left $\mathcal{C}$-approximation; see \cite{AS}.

\begin{lem}[{\cite[Theorem 1]{Hu}}]{\label{equivalence of relative n-torsionfree}} Let $\omega$ be an $R$-module with fadim$_R\omega\geq n+2$. For a module $M\in R$-${\rm mod}$, the following statements are equivalent.

$(1)$ $M$ is $\omega$-$n$-torsionfree.

$(2)$ There exists an exact sequence $\xymatrix@C=0.5cm{0 \ar[r] & M\ar[r]^{f_1} &{\omega_1}\ar[r]^{f_2~} &  \cdots \ar[r]^{f_n} &{\omega_n}}$ with $\omega_i\in{\rm add}\omega$, such that each ${\rm Im} f_i \rightarrow \omega_i$ is a left ${\rm add}\omega$-approximation of ${\rm Im} f_i$ for $1\leq i\leq n$.
\end{lem}

This description of $\omega$-$n$-torsionfree modules leads to the related notion of $\omega$-$n$-syzygy modules. A module $M$ is called an $\omega$-$n$-syzygy if there exists an exact sequence $0\rightarrow M\rightarrow X_{n-1}\rightarrow\cdots\rightarrow X_0\rightarrow N \rightarrow 0$ in $R$-${\rm mod}$ with $X_i$ in add$_R\omega$ for $0\leq i\leq n-1$. We denote by $\Omega_\omega^n(R)$ the full subcategory of $R$-${\rm mod}$ whose objects are  $\omega$-$n$-syzygy modules. It is clear that  $\mathcal{T}_\omega^n(R)\subseteq\Omega^n_\omega(R)$ by  Lemma \ref{equivalence of relative n-torsionfree}. But an $\omega$-$n$-syzygy module is not $\omega$-$n$-torsionfree in general.

 By the property of approximations, we have the following observation.
\begin{lem}{\label{lemma of relative torsionfree}} Let $\omega$ be an $R$-module with {\rm fadim}$_R\omega\geq n+2$ and ${\rm Ext}_R^i(\omega, \omega)=0$ for $1\leq i\leq n$ and $M$ be an $R$-module. Then $M$ is $\omega$-$n$-torsionfree if and only if there exists an exact sequence
$$\xymatrix@C=0.5cm{0 \ar[r] & M\ar[r]^{f_1} &{\omega_1}\ar[r]^{f_2~} &  \cdots \ar[r]^{f_n} &{\omega_n}\ar[r] & T_n \ar[r] & 0}$$
with $\omega_i\in{\rm add}_R\omega$ and  ${\rm Ext}^{j}_R(T_n, \omega)=0$ for $1\leq j\leq n$.
\end{lem}
\noindent{\bf Proof.} If $M$ is $\omega$-$n$-torsionfree, then there exists an exact sequence  $$\xymatrix@C=0.5cm{0 \ar[r] & M\ar[r]^{f_1} &{\omega_1}\ar[r]^{f_2~} &  \cdots \ar[r]^{f_n} &{\omega_n}\ar[r] & T_n \ar[r] & 0}$$
 with $\omega_i\in{\rm add}_R\omega$, such that ${\rm Im} f_i \rightarrow \omega_i$ is a left ${\rm add}_R\omega$-approximation of ${\rm Im} f_i$ for $1\leq i\leq n$ by Lemma \ref{equivalence of relative n-torsionfree}. Taking ${\rm Coker}f_i=T_i$, we have that ${\rm Ext}^1_R(T_i, \omega)=0$ for $1\leq i\leq n$. By dimension shifting,  ${\rm Ext}^j_R(T_i, \omega)\cong{\rm Ext}^{n-i+j}_R(T_n, \omega)$, and we obtain that ${\rm Ext}^{j}_R(T_n, \omega)=0$ for $1\leq j\leq n$.

Conversely, if there exists an exact sequence
$$\xymatrix@C=0.5cm{0 \ar[r] & M\ar[r]^{f_1} &{\omega_1}\ar[r]^{f_2~} & \cdots \ar[r]^{f_n} &{\omega_n}\ar[r] & T_n\ar[r] & 0}$$
with $\omega_i\in{\rm add}_R\omega$ and  ${\rm Ext}^{j}_R(T_n, \omega)=0$ for $1\leq j\leq n$, then ${\rm Ext}^1_R(T_i, \omega)\cong{\rm Ext}^{n-i+1}_R(T_n, \omega)=0$
for $1\leq i\leq n$ where $T_i={\rm Coker}f_i$. It follows that ${\rm Im} f_i \rightarrow \omega_i$ is a left ${\rm add}_R\omega$-approximation of ${\rm Im} f_i$ for $1\leq i\leq n$, and $M$ is $\omega$-$n$-torsionfree by Lemma \ref{equivalence of relative n-torsionfree} again.\hfill$\square$

A {\it ring extension} $S/R$ is a ring homomorphism $l: R\rightarrow S$ which is required to send the unit to the unit. Thus $R$ is the base ring and $S$ is the extension ring.  A ring extension is an \emph{algebra} if $R$ is commutative and $l$ factors $R\rightarrow C(S)\hookrightarrow S$, where $C(S)$ is the center of $S$. The natural bimodule structure of ${_RS_R}$ is given by $r\cdot s\cdot r^\prime:=l(r)\cdot s\cdot l(r^\prime)$. Similarly, one obtains the other natural module structures, for example $S_R$, ${_RS_S}$ and ${_SS_R}$, etc.

As a generalization of Frobenius algebras, Frobenius extensions were first introduced by Kasch  in \cite{Kas1}. A ring extension $S/R$ is said to be a \emph{Frobenius extension} if ${_RS}$ is finitely generated projective and ${_SS_R}\cong {\rm Hom}_R({_RS_S},{_RR_R})$ as $S$-$R$-bimodules. By \cite[Theorem 1.2]{Kad}, these conditions are equivalent to the corresponding properties on the opposite sides, that is, a ring extension $S/R$ is a Frobenius extension if and only if ${S_R}$ is finitely generated projective and ${_RS_S}\cong {\rm Hom}_{R^{{\rm op}}}({_SS_R},{_RR_R})$ as $R$-$S$-bimodules.

 Let $S/R$ be a ring extension and $M$ an $S$-module. Then $M$ is naturally an $R$-module, which is referred to as the \emph{underlying} $R$-module. There is a natural surjective map $\pi: S\otimes_RM\rightarrow M$ via $s\otimes m\mapsto sm $ for any $s\in S$ and $m\in M$. It is easy to check that $\pi$ is split as an $R$-homomorphism.  However, $\pi$ is not split as an $S$-homomorphism in general.

For a ring extension $S/R$, there is a {\it restriction functor} $Res: S$-${\rm Mod}\rightarrow R$-${\rm Mod}$ sending ${_SM}\mapsto {_RM}$, given by $r\cdot m:=l(r)\cdot m$. In the opposite direction, there are two natural functors as follows:

$(1)$ $\mathbb{T}={_SS}\otimes_R- :R$-${\rm Mod}\rightarrow S$-${\rm Mod}$ is given by ${_RM}\mapsto {_SS}\otimes_RM$.

$(2)$ $\mathbb{H}={\rm Hom}_R({_RS_S},-):R$-${\rm Mod}\rightarrow S$-${\rm Mod}$ is given by ${_RM}\mapsto {\rm Hom}_R({_RS_S},{_RM})$.

It is easy to check that both $(\mathbb{T},Res)$ and $(Res,\mathbb{H})$ are adjoint pairs.

By \cite[Theorem 1.2]{Kad} and \cite[Definition 2.1]{Kad}, a ring extension $S/R$ is a Frobenius extension if and only if the functors $\mathbb{T}$ and $\mathbb{H}$ above are naturally equivalent. Equivalently,  $(\mathbb{T}\cong\mathbb{H}, Res)$ is a Frobenius pair and $\mathbb{T}\cong\mathbb{H}$ are Frobenius functors \cite{CGN}.

A fundamental example of Frobenius extensions is the group
 algebra extension induced by a subgroup of  finite index. There are other examples of Frobenius extensions including Hopf subalgebras, finite extensions of
enveloping algebras of Lie super-algebras and finite extensions of enveloping algebras of Lie coloralgebras  etc \cite{FMS,Sch}. More examples of Frobenius extensions can be found in \cite[Example 2.4]{Z}.

Recall that an $R$-$R$-bimodule $M$ is \emph{centrally projective} over $R$ if ${_RM_R}\in{\rm add}({_RR_R})$; see \cite{H}.

\begin{exm} $(1)$ If $l: R\rightarrow S$ is a Frobenius extension and $S$ is an $R$-algebra, then  ${_RS_R}$ is centrally projective over $R$.

$(2)$ Let $G$ be a finite group and $H$ a subgroup of $G$.
For any ring $R$, the inclusion $R[H]\hookrightarrow R[G]$
is a Frobenius extension. Suppose that
$G=HC_G(H)$, where $C_G(H)$ is the centralizer of
$H$ in $G$. Choose coset representatives
$g_1,\ldots,g_m\in C_G(H)$, where $m=[G:H]$.
Since each $g_i$ commutes with the elements of
$R[H]$, we have an isomorphism of
$R[H]$-$R[H]$-bimodules
\[
{}_{R[H]}R[G]_{R[H]}
\cong
\bigoplus_{i=1}^{m}
{}_{R[H]}R[H]_{R[H]}.
\]
Thus $R[H]\hookrightarrow R[G]$ is a Frobenius
extension with ${}_{R[H]}R[G]_{R[H]}$ centrally
projective over $R[H]$; see also
\cite[Corollary~1.3]{BK}.

In particular, $R\hookrightarrow R[G]$ is a
Frobenius extension with ${}_RR[G]_R$
centrally projective over $R$.

$(3)$ Recall from \cite{P} that a ring extension $l:R\rightarrow S$ is called an {\em excellent extension} if
 \begin{enumerate}
	\item[(i)] $S$ is left $R$-projective, that is, if $_SM$ is a module and $_SN$ is a submodule then $_RN|_RM$ implies that $_SN|{_SM}$, where $N|M$ means that $N$ is a  direct summand of $M$.
	\item[(ii)] $S$ is a free normal extension of $R$, i.e., $S=\sum\limits_{i=1}^{n}s_iR$ and $S$ is free with common basis $\{s_1=1, s_2, \cdots, s_n\}$ as both a left $R$-module and a right $R$-module, such that $s_iR=Rs_i$ for $1\leq i\leq n$.
\end{enumerate}

By \cite[Lemma~4.7]{HS}, every excellent
extension is a Frobenius extension. Although
an excellent extension need not be centrally
projective in general, the following standard
examples satisfy this additional condition.

First, for any ring $R$, the canonical embedding
$R\hookrightarrow M_n(R)$ is an excellent
extension, and
\[
{}_RM_n(R)_R\cong({}_RR_R)^{n^2}.
\]

Second, let $K/k$ be a finite separable field
extension and $A$ a $k$-algebra. Then the scalar
extension
\[
A\longrightarrow A\otimes_k K
\]
is an excellent extension. Moreover,
\[
{}_A(A\otimes_k K)_A
\cong({}_AA_A)^{[K:k]}.
\]
Thus both extensions are Frobenius extensions
whose underlying bimodules are centrally
projective over the respective base rings.

\end{exm}

\section{Faithful dimensions and Wakamatsu tilting modules}

The faithful dimension was introduced first by Buan and Solberg in \cite{BS}. A Wakamatsu tilting module is just a self-orthogonal module with infinite faithful dimension.  We will investigate transfer properties of faithful dimensions under Frobenius extensions in this section.

\begin{lem}\label{transfer of ffdim} Let $S/R$  be a Frobenius extension with $_RS_R$ centrally projective over $R$. If {\rm fadim}$_R\omega\geq n+2$ and ${\rm Ext}_R^i(\omega, \omega)=0$ for $1\leq i\leq n$, then {\rm fadim}$_S({S}\otimes_R\omega)\geq n+2$.
\end{lem}

\noindent{\bf Proof.} By assumption, there exists an exact sequence
$$
0\longrightarrow R
\xrightarrow{f_1}\omega_1
\xrightarrow{f_2}\cdots
\xrightarrow{f_{n+2}}\omega_{n+2},
$$
where $\omega_i\in\operatorname{add}_R\omega$ and each inclusion
$\operatorname{Im}f_i\hookrightarrow\omega_i$ is a left
$\operatorname{add}_R\omega$-approximation for $1\leq i\leq n+2$. For each $i$, put $U_i=\operatorname{Im}f_i$ and
$T_i=\operatorname{Coker}f_i$. Then we have a short exact sequence
$$
0\longrightarrow U_i\longrightarrow\omega_i
\longrightarrow T_i\longrightarrow0.
$$
The approximation property implies that
$$
\operatorname{Hom}_R(\omega_i,\omega)
\longrightarrow
\operatorname{Hom}_R(U_i,\omega)
$$
is surjective. Since
$\operatorname{Ext}_R^1(\omega_i,\omega)=0$, the associated
long exact sequence yields
$\operatorname{Ext}_R^1(T_i,\omega)=0$.

Since ${}_RS_R\in\operatorname{add}({}_RR_R)$, we have
${}_R(S\otimes_R\omega)\in\operatorname{add}_R\omega$.
Consequently, by the adjoint isomorphism,
$$
\begin{aligned}
\operatorname{Ext}_S^1
(S\otimes_R T_i,S\otimes_R\omega)
&\cong
\operatorname{Ext}_R^1
(T_i,{}_R(S\otimes_R\omega))=0
\end{aligned}
$$
for $1\leq i\leq n+2$.

Applying the exact functor $S\otimes_R-$ to the short exact
sequence above, we obtain
$$
0\longrightarrow S\otimes_R U_i
\longrightarrow S\otimes_R\omega_i
\longrightarrow S\otimes_R T_i
\longrightarrow0.
$$
The preceding Ext vanishing shows that
$$
\operatorname{Hom}_S
(S\otimes_R\omega_i,S\otimes_R\omega)
\longrightarrow
\operatorname{Hom}_S
(S\otimes_R U_i,S\otimes_R\omega)
$$
is surjective. Since
$S\otimes_R\omega_i\in
\operatorname{add}_S(S\otimes_R\omega)$, the inclusion
$S\otimes_R U_i\longrightarrow S\otimes_R\omega_i$
is a left $\operatorname{add}_S(S\otimes_R\omega)$-approximation for $1\leq i\leq n+2$.

Finally, applying $S\otimes_R-$ to the original exact sequence
and using
\[
S\otimes_R R\cong S,
\qquad
S\otimes_R U_i\cong
\operatorname{Im}(S\otimes_R f_i),
\]
we obtain
$\operatorname{fadim}_S(S\otimes_R\omega)\geq n+2$. \hfill$\square$

 The following result shows that the Wakamatsu  tilting property is transferable under Frobenius extensions.

\begin{prop}{\label{WT transfer}} Let $S/R$  be a Frobenius extension with $_RS_R$ centrally projective over $R$. If $_R\omega$ is a Wakamatsu tilting module as an $R$-module, then so is $_S{S}\otimes_R\omega$ as an $S$-module.
\end{prop}
\noindent{\bf Proof.} Since $_R\omega$ is a Wakamatsu tilting module, we have {\rm fadim}$_R\omega=\infty$ and ${\rm Ext}_R^i(\omega,\omega)=0$ for any $i\geq 1$. Then {\rm fadim}$_S({S}\otimes_R\omega)=\infty$ by Lemma \ref{transfer of ffdim}. By assumption, ${_RS}\otimes_R\omega\in {\rm add}\omega$. It follows that
\begin{align*}{\rm Ext}_S^i(S\otimes_R\omega,S\otimes_R\omega)
&\cong {\rm Ext}_R^i(\omega,{\rm Hom}_S({_SS_R},S\otimes_R\omega))\\
&\cong {\rm Ext}_R^i(\omega,{_RS}\otimes_R\omega)=0
\end{align*}
for any $i\geq 1$ by adjoint isomorphisms.
\hfill$\square$

Let $S/R$  be a Frobenius extension with $_RS_R$ centrally projective over $R$ and $_R\omega$ a Wakamatsu tilting module. By Proposition \ref{WT transfer}, $_S{S}\otimes_R\omega$ is also a Wakamatsu tilting module. Then $_S({S}\otimes_R\omega)_\Gamma$ is a faithful balanced self-orthogonal bimodule, where $\Gamma={\rm End}_S({S}\otimes_R\omega)$. It is straightforward to  check that the map $\rho:T\rightarrow \Gamma$, defined by  $\rho(f): s\otimes w\mapsto s\otimes (w\cdot f)$ for any $f\in T$, $s\in S$ and $w\in\omega$, is a ring homomorphism. We will show that $\rho:T\rightarrow \Gamma$ is also a Frobenius extension provided that pd$(\omega_T)$ is finite.

\begin{lem}{\label{generalized tilting and injective}} Let $_R\omega$ be a Wakamatsu tilting module and $T={\rm End}_R(\omega)$. For any injective $R$-module $I$, we have ${\rm Tor}_i^T(\omega_T, {\rm Hom}_R(\omega_T, I))=0$ for any $i\geq 1$.
\end{lem}
\noindent{\bf Proof.} Since $\omega_T$ is finitely generated and $I$ is injective, the canonical map
 $\varphi: {_R\omega}\otimes_T{\rm Hom}_R(_R\omega_T, I)\rightarrow {\rm Hom}_R({\rm Hom}_{T^{\rm op}}({_R\omega},{_R\omega}), I)$ given by
 $w\otimes f\mapsto \varphi(w\otimes f)(g)=(f\cdot g)(w)$
 is an isomorphism; see \cite[Page12]{Chris}. Since $_R\omega$ is a Wakamatsu tilting module, we have ${\rm Hom}_{T^{\rm op}}({_R\omega},{_R\omega})\cong {_RR_R}$. And therefore $_R\omega\otimes_T{\rm Hom}_R(_R\omega_T, I)\cong {_RI}$.

 Take a projective resolution of $\omega_T$ in mod-$T$ as follows
 \begin{equation}{\label{3.1}}
\xymatrix@C=0.5cm{\cdots \ar[r] & P_1\ar[r] &P_0\ar[r] & {_R\omega_T}  \ar[r] & 0.}
\end{equation}
Since $\omega_T$ is self-orthogonal, we have
 the following exact sequence $$0\rightarrow {\rm Hom}_{T^{\rm op}}(_R\omega_T,{_R\omega})\cong {_RR_R}\rightarrow {\rm Hom}_{T^{\rm op}}(P_0, {_R\omega_T})\rightarrow {\rm Hom}_{T^{\rm op}}(P_1, {_R\omega_T})\rightarrow\cdots$$
 in $R$-mod after applying the functor ${\rm Hom}_{T^{\rm op}}(-,{_R\omega})$ to the exact sequence (\ref{3.1}).
 Applying the exact functor ${\rm Hom}_R(-, I)$ to the sequence above, we get the following exact sequence
\begin{equation}{\label{3.2}}
\xymatrix@C=0.3cm{\cdots \ar[r] &{\rm Hom}_R({\rm Hom}_{T^{\rm op}}({P_1},{_R\omega_T}), I) \ar[r] &{\rm Hom}_R({\rm Hom}_{T^{\rm op}}({P_0},{_R\omega_T}), I)\ar[r] & {\rm Hom}_R(R, I)\ar[r] & 0.}
\end{equation}

On the other hand, applying  the functor $-\otimes_T{\rm Hom}_R({_R\omega_T, I})$ to the exact sequence (\ref{3.1}), we get the following complex
\begin{equation}{\label{3.3}}
\xymatrix@C=0.3cm{\cdots \ar[r] & P_1\otimes_T{\rm Hom}_R({_R\omega_T, I})\ar[r] & P_0\otimes_T{\rm Hom}_R({_R\omega_T, I})\ar[r] & \omega\otimes_T{\rm Hom}_R({_R\omega_T, I})\ar[r] & 0.}
\end{equation}
Comparing the exact sequence (\ref{3.2}) with the complex (\ref{3.3}), we see that the sequence  (\ref{3.3}) is also exact since $P_i\otimes_T{\rm Hom}_R({_R\omega_T, I})\cong {\rm Hom}_R({\rm Hom}_{T^{\rm op}}({P_i},{_R\omega_T}), I)$ for any $i\geq 0$ and these isomorphisms are natural. Thus, we get that $${\rm Tor}_i^T(\omega_T, {\rm Hom}_R(\omega_T, I))=0$$ for $i\geq 1$.
\hfill$\square$

\begin{prop} Let $S/R$  be a Frobenius extension with $_RS_R$ centrally projective over $R$ and $_R\omega$ be a Wakamatsu tilting module with $T={\rm End}_R(\omega)$. If ${\rm pd}(\omega_T)$ is finite, then $_R\omega\otimes_T{\rm Hom}_R(_R\omega_T, S\otimes_R\omega)\cong {_RS\otimes_R\omega}$.
\end{prop}
\noindent{\bf Proof.}
Since $_RS_R\in{\rm add}({_RR_R})$ by assumption, we have $_RS\otimes_R\omega\in{\rm add}_R\omega$.  Since  ${\rm Ext}_R^i({_R\omega},_R\omega )=0$ for $i\geq 1$, we obtain that ${\rm Ext}_R^i(\omega,_RS\otimes_R\omega )=0$ for $i\geq 1$.

 By Lemma \ref{generalized tilting and injective}, we know that the bimodule $_R\omega_T$ and the $R$-module $_RS\otimes_R\omega$  satisfy the hypotheses of Theorem 1.14 in \cite{Mi}.
Therefore, we have
$$_R\omega\otimes_T{\rm Hom}_R(_R\omega_T, S\otimes_R\omega)\cong {_RS\otimes_R\omega}$$
as a direct consequence of \cite[Theorem 1.14]{Mi}.
\hfill$\square$

The following is the main result of this section.

\begin{thm}{\label{The main thm A}} Let $S/R$  be a Frobenius extension with $_RS_R$ centrally projective over $R$ and $_R\omega$ be a Wakamatsu tilting module. Set $T={\rm End}_R(\omega)$ and
$\Gamma={\rm End}_S({S}\otimes_R\omega)$. If ${\rm pd}(\omega_T)$ is finite, then the ring homomorphism $\rho: T\rightarrow \Gamma$ is a Frobenius extension.
\end{thm}
\noindent{\bf Proof.} As a $T$-module, we have
\begin{align*}_T\Gamma
&={\rm Hom}_S(S\otimes_R\omega_T,S\otimes_R\omega)\\
&\cong{\rm Hom}_R(\omega_T, {\rm Hom}_S({_SS_R}, S\otimes_R\omega))\\
&\cong{_T({\rm Hom}_R(\omega,{_RS}\otimes_R\omega)}).
\end{align*}
Since $_RS_R\in {\rm add}(_RR_R)$ by assumption, we have ${_RS}\otimes_R\omega\in{\rm add}(_R\omega)$. It follows that $_T\Gamma$ is a finitely generated projective $T$-module.

On the other hand, one has 
\begin{align*}{\rm Hom}_T({_T\Gamma_{\Gamma}},{_TT_T})
&={\rm Hom}_T({\rm Hom}_S(S\otimes_R\omega_T,S\otimes_R\omega),{\rm Hom}_R(\omega_T,\omega_T))\\
&\stackrel{\rm adj. iso.}{\cong}{\rm Hom}_T({\rm Hom}_R(\omega_T,S\otimes_R\omega),{\rm Hom}_R(\omega_T,\omega_T))\\
&\stackrel{\rm adj. iso.}{\cong}{\rm Hom}_R(\omega\otimes_T{\rm Hom}_R(\omega_T,S\otimes_R\omega), \omega_T)\\
&\stackrel{\rm Prop. 3.4}{\cong}{\rm Hom}_R(S\otimes_R\omega,{_R\omega_T})\\
&\quad\cong{\rm Hom}_R({_RS\otimes_SS\otimes_R\omega},{_R\omega_T})\\
&\stackrel{\rm adj. iso.}{\cong}{\rm Hom}_S(S\otimes_R\omega,{\rm Hom}_R({_RS_S},{_R\omega_T}))\\
&\quad\stackrel{\rm F.E}{\cong}{\rm Hom}_S(S\otimes_R\omega,S\otimes_R\omega_T)\\
&\quad\cong{_\Gamma\Gamma_T},
\end{align*}
where ``adj. iso." means ``adjoint isomorphisms". Here the third isomorphism is induced by the canonical
evaluation isomorphism
$\omega\otimes_T\Hom_R(\omega,S\otimes_R\omega)\xrightarrow{\sim}S\otimes_R\omega$
from Proposition~3.4, while the sixth
isomorphism is given by the Frobenius
isomorphism
$\Hom_R({}_RS_S,-)\cong S\otimes_R-$.
All these isomorphisms are natural and compatible with the
left $\Gamma$-action and the right $T$-action. Hence
\[
{}_\Gamma\Hom_T({}_T\Gamma,{}_TT_T)
\cong{}_\Gamma\Gamma_T
\]
as $\Gamma$-$T$-bimodules.

Since $_T\Gamma$ is  finitely generated projective, $\rho: T\rightarrow \Gamma$ is a Frobenius extension by \cite[Theorem 1.2]{Kad}.
\hfill$\square$

\section{$\omega$-$n$-Torsionfree Modules}

In this section, we assume throughout that fadim$_R\omega\geq n+2$ and ${\rm Ext}_R^i(\omega, \omega)=0$ for $1\leq i\leq n$. Then $_R\omega_T$ is faithful and balanced,
and ${\rm Ext}_R^i(\omega, \omega)=0={\rm Ext}_{T^{\rm op}}^i(\omega,\omega)$ for $1\leq i\leq n$, where $T={\rm End}_R(\omega)$. We will show that the relative $n$-torsionfreeness is preserved under Frobenius extensions.

The following result is a relative version of \cite[Proposition 2.3]{Z1}, which gives a recurrence relation of $\omega$-$n$-torsionfreeness.
\begin{prop}\label{recurrence of k-torsionfree} Let $\xymatrix@C=0.5cm{0\ar[r] & M\ar[r]^{f} &X\ar[r] & N\ar[r] &0}$ be a short exact sequence in $R$-${\rm mod}$ with $X\in{\rm add}\omega$ and $f^\omega$ being epic and let $n\geq 1$ be an integer.  Then $M$ is $\omega$-$n$-torsionfree if and only if $N$ is $\omega$-$(n-1)$-torsionfree.
\end{prop}
\noindent{\bf Proof.} $(\Rightarrow)$ If $n=1$, it is trivial.

For the case of  $n=2$, we will show that if $M$ is $\omega$-reflexive then $N$ is $\omega$-torsionless. Since $f^\omega$ is epimorphic, applying the functor ${\rm Hom}_R(-, \omega)=(-)^\omega$, we get an exact sequence
\begin{equation}{\label{4.1}} \xymatrix@C=0.5cm{0\ar[r] & N^\omega\ar[r] &X^\omega\ar[r]^{f^\omega} & M^\omega\ar[r] &0.}
\end{equation}

Consider the following commutative diagram with exact rows
$$\xymatrix{
 0\ar[r]& M\ar[d]_{\sigma_M} \ar[r]^{f} & X\ar[d]_{\sigma_X} \ar[r] &N\ar[d]_{\sigma_N}\ar[r]& 0\\
  0\ar[r]& M^{\omega\omega} \ar[r]^{f^{\omega\omega}} &  X^{\omega\omega}\ar[r] & N^{\omega\omega}.}$$
By the Snake Lemma,  $\sigma_N$ is a monomorphism since $\sigma_M$ and $\sigma_X$ are isomorphisms. Hence $N$ is an $\omega$-torsionless module.

Now assume that $n\geq3$. Since $M$ is
$\omega$-$n$-torsionfree, it is $\omega$-reflexive.
Moreover, by dimension shifting and the standing assumptions,
we have
\[
\operatorname{Ext}_{T^{\mathrm{op}}}^{j}(M^\omega,\omega)
\cong
\operatorname{Ext}_{T^{\mathrm{op}}}^{j+2}
(\operatorname{Tr}_\omega M,\omega)
=0
\]
for $1\leq j\leq n-2$. In particular,
$\operatorname{Ext}_{T^{\mathrm{op}}}^{1}(M^\omega,\omega)=0$.

Applying $\operatorname{Hom}_{T^{\mathrm{op}}}(-,\omega)$
to (4.1), we therefore obtain a short exact sequence
\[
0\longrightarrow M^{\omega\omega}
\xrightarrow{f^{\omega\omega}}X^{\omega\omega}
\longrightarrow N^{\omega\omega}
\longrightarrow0.
\]
Since $X\in\operatorname{add}_R\omega$, the evaluation map
$\sigma_X$ is an isomorphism. As $\sigma_M$ is also an
isomorphism, comparison with the original short exact sequence
shows that $\sigma_N$ is an isomorphism.
Hence $N$ is $\omega$-reflexive, and thus
$\operatorname{Ext}_{T^{\mathrm{op}}}^{i}
(\operatorname{Tr}_\omega N,\omega)=0$ for $i=1,2$. This proves the assertion when $n=3$.

Suppose now that $n\geq4$.
Since $X^\omega$ is a projective right $T$-module,
the exact sequence (4.1), together with dimension shifting,
gives
\[
\begin{aligned}
\operatorname{Ext}_{T^{\mathrm{op}}}^{j+2}
(\operatorname{Tr}_\omega N,\omega)
&\cong
\operatorname{Ext}_{T^{\mathrm{op}}}^{j}(N^\omega,\omega)\\
&\cong
\operatorname{Ext}_{T^{\mathrm{op}}}^{j+1}(M^\omega,\omega)\\
&=0
\end{aligned}
\]
for $1\leq j\leq n-3$.
Consequently,
$\operatorname{Ext}_{T^{\mathrm{op}}}^{i}
(\operatorname{Tr}_\omega N,\omega)=0$
for $1\leq i\leq n-1$.
Therefore, $N$ is $\omega$-$(n-1)$-torsionfree.

$(\Leftarrow)$ If $N$ is $\omega$-$(n-1)$-torsionfree, then there exists an exact sequence $$\xymatrix@C=0.5cm{0 \ar[r] & N\ar[r]^{f_1}& {X_1}\ar[r]^{f_2~} &  \cdots \ar[r]^{f_{n-1}} &{X_{n-1}}}$$ with each $X_i\in{\rm add}\omega$, such that ${\rm Im} f_i \hookrightarrow X_i$ is a left ${\rm add}\omega$-approximation of ${\rm Im} f_i$ for $1\leq i\leq n-1$ by Lemma \ref{equivalence of relative n-torsionfree}. Consider the exact sequence $$\xymatrix@C=0.5cm{0 \ar[r] & M \ar[r]^f & X\ar[r]^{f_1} &{X_1}\ar[r]^{f_2~} &  \cdots \ar[r]^{f_{n-1}} &{X_{n-1}},}$$ where $M(\cong{\rm Im}f)\rightarrow X$ is a left ${\rm add}\omega$-approximation (since $f^\omega$ is an epimorphism) and ${\rm Im} f_i \hookrightarrow X_i$ is a left ${\rm add}\omega$-approximation of ${\rm Im} f_i$ for $1\leq i\leq n-1$. Hence $M$ is $\omega$-$n$-torsionfree by Lemma \ref{equivalence of relative n-torsionfree} again.  \hfill$\square$

\begin{lem}\label{lemma of add} Let $S/R$ be a Frobenius extension with $_RS_R$ centrally projective over $R$ and $M$ be an $S$-module. If $M\in{\rm add}_S({S\otimes_R\omega})$, then $_RM\cong {_RS}\otimes_SM\in{\rm add}{_R\omega}$.
\end{lem}
\noindent{\bf Proof.} If $_SM\in{\rm add}_S({S}\otimes_R\omega)$, then we have
$_RM\in{\rm add}_R({S}\otimes_R\omega)$. Since $_RS_R$ is centrally projective over $R$, we have $_RS\otimes_R\omega\in {\rm add}{_R\omega}$. Hence $_RM\in{\rm add}{_R\omega}$. \hfill$\square$

\begin{prop}{\label{equivalence I}} Let $S/R$  be a Frobenius extension with $_RS_R$ centrally projective over $R$ and $M$ be an $R$-module. If $M$ is $\omega$-$n$-torsionfree as an $R$-module, then $S\otimes_RM$ is also $({S\otimes_R\omega})$-$n$-torsionfree as an $S$-module.
\end{prop}
\noindent{\bf Proof.} Suppose that ${_RM}$ is an $\omega$-$n$-torsionfree $R$-module. Then there exists an exact sequence
$$\xymatrix@C=0.5cm{0 \ar[r] & {_RM}\ar[r]^{f_1} &{\omega_1}\ar[r]^{f_2~} &{\omega_2}\ar[r] &\cdots \ar[r]^{f_n} &{\omega_n}\ar[r]& {T_n}\ar[r]& 0}$$
in $R$-${\rm mod}$ with $\omega_i\in{\rm add}_R\omega$ and ${\rm Ext}^j_R(T_n, \omega)=0$ for $1\leq j\leq n$ by Lemma \ref{lemma of relative torsionfree}.
Since $S_R$ is finitely generated projective by assumption, the functor $S\otimes_R-$ is exact. We get the following exact sequence
$$\xymatrix@C=0.5cm{0 \ar[r] & {S\otimes_RM}\ar[r]^{S\otimes_Rf_1} &{S\otimes_R\omega_1}\ar[r]^{S\otimes_Rf_2} & {S\otimes_R\omega_2}\ar[r] & \cdots \ar[r]^{S\otimes_Rf_n} &{S\otimes_R\omega_n}\ar[r]& {S\otimes_RT_n}\ar[r]& 0}$$
in $S$-${\rm mod}$ with $S\otimes_R\omega_i\in{\rm add}_S(S\otimes_R\omega)$ for $1\leq i\leq n$ and ${S\otimes_RT_n}\cong {\rm Coker}(S\otimes_Rf_n)$.
By the adjoint isomorphism, for any $i\geq 0$, we have
\begin{align*}{\rm Ext}^i_S(S\otimes_RT_n,{_SS\otimes_R\omega})
&\cong{\rm Ext}^{i}_R({T_n}, {\rm Hom}_S({_SS_R}, S\otimes_R\omega))\\
&\cong {\rm Ext}^i_R({T_n}, {_RS\otimes_R\omega}).
\end{align*}
Since $_RS\otimes_R\omega\in {\rm add}_R\omega$ by  assumption, we have $$0={\rm Ext}^i_R({T_n}, {S\otimes_R\omega})\cong{\rm Ext}^i_S(S\otimes_RT_n,S\otimes_R\omega)$$ for $1\leq i\leq n$.
Then $S\otimes_RM$ is $(S\otimes_R\omega)$-$n$-torsionfree as an $S$-module by Lemma \ref{lemma of relative torsionfree}.
\hfill$\square$

The following is the main result of this section.

\begin{thm}{\label{equivalence II}} Let $S/R$  be a Frobenius extension with $_RS_R$ centrally projective over $R$ and $M$ be an $S$-module. Then $M$ is $ (S\otimes_R\omega)$-$n$-torsionfree as an $S$-module if and only if $M$ is $\omega$-$n$-torsionfree as the underlying $R$-module.
\end{thm}
\noindent{\bf Proof.} $(\Rightarrow)$ Suppose that $M$ is an $(S\otimes_R\omega)$-$n$-torsionfree $S$-module. By Lemma \ref{lemma of relative torsionfree}, there exists an exact sequence
\begin{equation}{\label{4.2}}
\xymatrix@C=0.5cm{0 \ar[r] & {_SM}\ar[r]^{f_1} &{X_1}\ar[r]^{f_2~} &  \cdots \ar[r]^{f_n} &{X_n}\ar[r] & T_n \ar[r] & 0}
\end{equation}
in $S$-${\rm mod}$ with $X_i\in{\rm add}{_S(S\otimes_R\omega)}$ such that ${\rm Ext}^{j}_S(T_n, S\otimes_R\omega)=0$ for $1\leq j\leq n$, where $T_n={\rm Coker}f_n$.
 Applying the restriction functor $_RS\otimes_S-$ to (\ref{4.2}), we get the following exact sequence
\begin{equation}
\xymatrix@C=0.5cm{0 \ar[r] & {_RM}\ar[r]^{f_1} &{X_1}\ar[r]^{f_2~} &  \cdots \ar[r]^{f_n} &{X_n}\ar[r] & T_n \ar[r] & 0}
\end{equation}
in $R$-${\rm mod}$. By Lemma \ref{lemma of add}, we have $X_i\in{\rm add}_R\omega$ for $1\leq i\leq n$. Since $S/R$  is a Frobenius extension,  we have ${\rm Hom}_R({_RS_S},-)\cong{_SS\otimes_R-}$. By the adjoint isomorphism, for $1\leq j\leq n$, we get
\begin{align*}{\rm Ext}_R^j(T_n,\omega)
&\cong {\rm Ext}_R^j({_RS\otimes_ST_n},\omega)\\
&\cong {\rm Ext}_S^j(T_n,{\rm Hom}_R({_RS_S},\omega))\\
&\cong {\rm Ext}_S^j(T_n,{_SS\otimes_R\omega})\\
&=0.
\end{align*}
Thus, $M$ is an $\omega$-$n$-torsionfree $R$-module by Lemma \ref{lemma of relative torsionfree}.

$(\Leftarrow)$ For the case of  $n=1$, $M$ is $\omega$-torsionless, which is equivalent to saying that $M$ is cogenerated by $_R\omega$, that is, there exists an exact sequence
$\xymatrix@C=0.5cm{0 \ar[r] & M\ar[r]^{f} &\omega^{I}}$ for some index set $I$. By assumption, $_SS\otimes_R-$ is an exact functor, we have  an exact sequence
 $\xymatrix@C=0.5cm{0 \xrightarrow{}{_SS\otimes_RM\cong{\rm Hom}_R({_RS_S},M)}\xrightarrow{S\otimes_Rf}{S\otimes_R\omega^I}}$ in $S$-mod. Since ${S\otimes_R-}\cong{\rm Hom}_R({_RS_S},-)$, we obtain
\begin{align*}
{S\otimes_R\omega^I}
&\cong{\rm Hom}_R({_RS_S},\omega^I)\\
&\cong\prod\limits_{i\in I}{\rm Hom}_R({_RS_S},\omega)\\
&\cong\prod\limits_{i\in I}{S\otimes_R\omega}\\
&\cong({S\otimes_R\omega})^I.
\end{align*}

 Note that there is an $S$-monomorphism $i:{_SM}\rightarrow {\rm Hom}_R({_RS_S},{_RM})$ via $i(m)(s)=sm$ for any $m\in M$ and $s\in S$. Since $S\otimes_RM\cong{\rm Hom}_R({_RS_S},M)$, we have a monomorphism, which we still denote by $i$, from $M$ to $S\otimes_RM$. It follows that there is  an exact sequence
$0 \rightarrow{_SM}\xrightarrow{(S\otimes_Rf) \cdot i}({S\otimes_R\omega})^I$.  Hence ${_SM}$ is cogenerated by $S\otimes_R\omega$, and so $M$ is $S\otimes_R\omega$-torsionless as an $S$-module.

Now we assume that $n>1$. Since $M$ is an $\omega$-$n$-torsionfree module as an $R$-module, there exists an exact sequence
$\xymatrix@C=0.5cm{0 \xrightarrow{}{_RM}\xrightarrow{f_1}{\omega_1}\xrightarrow{f_2~}{\omega_2}\xrightarrow{}\cdots}\xrightarrow{f_n}{\omega_n}$
in $R$-mod with $\omega_i\in{\rm add}_R\omega$ and ${\rm Im}f_i\hookrightarrow \omega_i$ being a left ${\rm add}_R\omega$-approximation of ${\rm Im}f_i$ for $1\leq i\leq n$ by Lemma \ref{equivalence of relative n-torsionfree}.
Applying  the exact functor $_SS\otimes_R-$, we get the following exact sequence
$$\xymatrix@C=0.8cm{0 \xrightarrow{}{_SS\otimes_RM}\xrightarrow{S\otimes_Rf_1}{S\otimes_R\omega_1}\xrightarrow{S\otimes_Rf_2}{S\otimes_R\omega_2}\xrightarrow{} \cdots\xrightarrow{S\otimes_Rf_n}{S\otimes_R\omega_n}}$$
in $S$-mod with $S\otimes_R\omega_i\in {\rm add}_S(S\otimes_R\omega)$  for $1\leq i\leq n$.
By Proposition \ref{equivalence I} and its proof, we have that $S\otimes_RM$ is $(S\otimes_R\omega)$-$n$-torsionfree as an $S$-module and ${\rm Im}(S\otimes_Rf_i)\hookrightarrow {S\otimes_R\omega_i}$ is a left ${\rm add}_S(S\otimes_R\omega)$-approximation for $1\leq i\leq n$.

Note that there is an $S$-epimorphism $\pi:{S\otimes_RM}\rightarrow {_SM}$ via $\pi(s\otimes m)=sm$ for any $m\in M$ and $s\in S$, which is split  as an $R$-homomorphism. Meanwhile, there is an $S$-monomorphism $i:{_SM}\rightarrow {\rm Hom}_R({_RS_S},{_RM})$ via $i(m)(s)=sm$ for any $m\in M$ and $s\in S$.
Since ${_SS\otimes_RM}\cong {\rm Hom}_R({_RS_S},M)$, there is an $S$-monomorphism, which we still denote  by  $i:{_SM}\hookrightarrow {S\otimes_RM}$ and it is split as a homomorphism of $R$-modules. It follows that there exist an $R$-homomorphism $p:{S\otimes_RM}\rightarrow M$ such that $p\cdot i=1_M$.
Set $f=S\otimes_Rf_1$. Then there exists an exact sequence
\begin{equation}{\label{4.4}}
\xymatrix@C=0.5cm{0 \xrightarrow{} {_SM}\xrightarrow{f\cdot i}{S\otimes_R\omega_1}\xrightarrow{}T\xrightarrow{}0,}
\end{equation}
where $T\cong {\rm Coker}(f\cdot i).$

  Applying the restriction functor to (\ref{4.4}), we get an exact sequence
  $$\xymatrix@C=0.5cm{0 \xrightarrow{} {_RM}\xrightarrow{f\cdot i}{S\otimes_R\omega_1}\xrightarrow{}T\xrightarrow{}0}$$ in $R$-mod with ${_RS\otimes_R\omega_1}\in {\rm add}_R\omega$. Let $h:M\rightarrow \omega$ be an $R$-homomorphism. Then $h\cdot p$ is an $R$-homomorphism from ${_RS\otimes_RM}$ to $\omega$, where $p$ is the split $R$-epimorphism from ${_RS\otimes_RM}$ to $M$. It follows that ${_RS\otimes_RM}$ is $\omega$-$n$-torsionfree and $f$ (as an $R$-homomorphism) is a left ${\rm add}(_R\omega)$-approximation of ${_RS\otimes_RM}$. Hence,
 there exists an $R$-homomorphism $g:{S\otimes_R\omega_1}\rightarrow \omega$ such that $g\cdot f=h\cdot p$
 $$\xymatrix{
  & _R\omega  &  \\
  0\ar[r]& {_RM} \ar[u]_{\forall h}\ar[r]^{f\cdot i} &  {S\otimes_R\omega_1}\ar@{..>}[ul]_{\exists g} \ar[r] & {_RT} \ar[r] & 0 \\
  &{_RS\otimes_RM} \ar[u]_{p} \ar[ur]_{f} & . }$$
Then, $g\cdot(f\cdot i)=(g\cdot f)\cdot i=(h\cdot p)\cdot i=h\cdot(p\cdot i)=h$. It follows that $$f\cdot i:{_RM}\hookrightarrow{_RS\otimes_R\omega_1}$$ is a left ${\rm add}_R\omega$-approximation of $_RM$.
Since $_RM$ is $\omega$-$n$-torsionfree, it follows that $T$ is $\omega$-$(n-1)$-torsionfree as an $R$-module by Proposition \ref{recurrence of k-torsionfree}. By the induction hypothesis, $T$ is also $(S\otimes_R\omega)$-$(n-1)$-torsionfree as an $S$-module. Since $f\cdot i:{_RM}\hookrightarrow{_RS\otimes_R\omega_1}$ is a left ${\rm add}_R\omega$-approximation of $M$, we have
\begin{align*}0
&={\rm Ext}^1_R({_RT},\omega)\\
&\cong{\rm Ext}^{1}_R({_RS\otimes_ST},\omega)\\
&\cong{\rm Ext}^{1}_S({_ST}, {\rm Hom}_R({_RS_S}, \omega))\\
&\cong {\rm Ext}^1_S({_ST}, {S\otimes_R\omega}).
\end{align*}
So we get that ${\rm Hom}_S({_SS\otimes_R\omega_1}, {S\otimes_R\omega})\rightarrow{\rm Hom}_S({_SM},{S\otimes_R\omega})\rightarrow 0$ is exact after applying the functor Hom$_S(-, {_SS\otimes_R\omega})$ to the exact sequence (\ref{4.4}).
Hence $f\cdot i:{_SM}\hookrightarrow{_SS\otimes_R\omega_1}$ is a left ${\rm add}_S(S\otimes_R\omega)$-approximation of ${_SM}$. Since $T$ is  an $({S\otimes_R\omega})$-$(n-1)$-torsionfree $S$-module, we conclude that ${_SM}$ is an $({S\otimes_R\omega})$-$n$-torsionfree module by Proposition \ref{recurrence of k-torsionfree}. This completes the proof. \hfill$\square$

Recall that a subcategory $\mathcal{X}$ of $R$-mod is called \emph{extension-closed}, if the middle term $X$ of any short exact sequence $0\rightarrow X'\rightarrow X\rightarrow X''\rightarrow 0$ is in $\mathcal{X}$, provided the end terms $X', X''$ are in $\mathcal{X}$. Extension-closed subcategories are of particular interest,  especially when they are  covariantly finite. We will investigate the transfer of the extension-closedness for the subcategory of $\omega$-$n$-torsionfree modules.

\begin{prop}{\label{equivalent of extension closure}}Let $S/R$  be a Frobenius extension with $_RS_R$ centrally projective over $R$ and $n$ be a positive integer. If  $\mathcal{T}^n_\omega(R)$ is extension-closed, then so is $\mathcal{T}^n_{S\otimes_R\omega}(S)$.
\end{prop}
\noindent{\bf Proof.} Let $0\rightarrow A\rightarrow B\rightarrow C\rightarrow 0$ be an exact sequence in $S$-mod such that $A, C$ are $S\otimes_R\omega$-$n$-torsionfree modules. Applying  the restriction  functor $_RS\otimes_S-$, we have the  exact sequence
$0\rightarrow A\rightarrow B\rightarrow C\rightarrow 0$
in $R$-mod. It follows from Theorem \ref{equivalence II} that $A$ and $C$ are $\omega$-$n$-torsionfree as $R$-modules.
By assumption, $B$ is also $\omega$-$n$-torsionfree as an $R$-module. So $B$ is $S\otimes_R\omega$-$n$-torsionfree as an $S$-module by Theorem \ref{equivalence II} again.   \hfill$\square$

\begin{cor} Let $S/R$  be a Frobenius extension with $_RS_R$ centrally projective over $R$ and $n$ a positive integer. If $\Omega^i_\omega(R)$ is extension-closed for $1\leq i\leq n$, then so is $\Omega^i_{S\otimes_R\omega}(S)$ for $1\leq i\leq n$. In this case, we have $\mathcal{T}_{\omega}^i(R)=\Omega^i_\omega(R)$ and $\mathcal{T}_{S\otimes_R\omega}^i(S)=\Omega^i_{S\otimes_R\omega}(S)$ for any $1\leq i\leq n$.
\end{cor}
\noindent{\bf Proof.} By Proposition \ref{equivalent of extension closure} and Theorem 3.6 in \cite{Hu1}, the assertion follows.
\hfill$\square$

\section{Modules of generalized G-dimension zero related to $\omega$}

In this section, $_R\omega$ is always a Wakamatsu tilting module, equivalently, $_R\omega_T$ is a faithful balanced self-orthogonal bimodule, where $T={\rm End}(_R\omega)$ is the endomorphism ring of $_R\omega$. As a generalization of modules with G-dimension zero, Auslander and Reiten introduced the notion of modules with generalized G-dimension zero with respect to $\omega$ in \cite{AR91}.
 A module $M$ in $R$-mod is said to have \emph{generalized Gorenstein dimension zero} with respect to $\omega$, denoted by G-dim$_\omega(M)=0$, if the following two conditions hold:
(1) $M$ is $\omega$-reflexive;
(2) ${\rm Ext}^i_R(M, \omega)=0={\rm Ext}^i_{T^{\rm op}}(M^\omega, \omega)$ for any $i\geq 1$. We denote by $\mathcal{G}_\omega$ the subcategory of $R$-${\rm mod}$ consisting of all modules of generalized G-dimension zero.

If $_R\omega_T={_RR_R}$, the module of generalized G-dimension zero is just the module of G-dimension zero. We mention that the notion of the module of G-dimension zero is the origin of Gorenstein homological algebra.
We denote the subcategory of $R$-${\rm mod}$ consisting of all modules $M$ with ${\rm Ext}^i_R(M, \omega)=0$ for all $i>0$ by ${^\bot_R\omega}$, and it is called \emph{the left orthogonal class} of $_R\omega$. An $R$-module $M$  is called  \emph{$\omega$-$\infty$-torsionfree} if it is $\omega$-$n$-torsionfree for any positive integer $n$. By definition of modules with generalized G-dimension zero, we have the following simple observation.

 \begin{lem}{\label{lem G G-dimension zero}}
 Let $M$ be an $R$-module. Then {\rm G}-${\rm dim}_\omega(M)=0$ if and only if  $M\in{^\bot_R\omega}$ and $M$ is $\omega$-$\infty$-torsionfree.
 \end{lem}
 \noindent{\bf Proof.} Let
 $\xymatrix@C=0.5cm{P_1\ar[r]^{f} & P_0\ar[r] & M  \ar[r] & 0}$ be a projective resolution of $_RM$ in $R$-${\rm mod}$.  Applying the functor $(-)^\omega={\rm Hom}_R(-,\omega)$, we have an exact sequence $$\xymatrix@C=0.5cm{0 \ar[r] & {M^\omega}\ar[r] &P_0^\omega\ar[r]^{f^\omega}&P_1^\omega \ar[r] &{\rm Tr}_\omega M\ar[r] & 0}.$$ By dimension shifting, we have ${\rm Ext}^i_{T^{\rm op}}(M^\omega, \omega)\cong{\rm Ext}^{i+2}_{T^{\rm op}}({\rm Tr}_\omega M, \omega)$.
 Then the conditions $M$ is $\omega$-reflexive and ${\rm Ext}^i_{T^{\rm op}}(M^\omega, \omega)=0$ for all $i\geq1$ are equivalent to  $M$ being $\omega$-$\infty$-torsionfree.
 \hfill$\square$

 A central notion in  Gorenstein homological algebra is that of a Gorenstein  projective module, introduced by Enochs and Jenda in \cite{EJ} as a generalization of modules of G-dimension zero. There are several different terminologies in the literature for these modules, such as totally reflexive modules, maximal Cohen-Macaulay modules. Modules of generalized G-dimension zero provide a  relative version of these modules and  also play an important role in relative homological algebra.

Chen proved  the total reflexivity of modules is preserved under the totally reflexive extension in \cite{CH}, where the totally reflexive extension is a generalization of Frobenius extensions.  Ren \cite{R} and Zhao \cite{Z} respectively proved that the Gorenstein projectivity of modules is preserved under Frobenius extensions by constructing  projective coresolutions.

In this section, we will show that a module over the extension ring is of generalized G-dimension zero if and only if its underlying module over the base ring has generalized G-dimension zero under Frobenius extensions.

 First, we record two immediate consequences of  Proposition \ref{equivalence I} and Theorem \ref{equivalence II}.
\begin{cor}{\label{torsionfree equ I}} Let $S/R$  be a Frobenius extension with $_RS_R$ centrally projective over $R$ and $M$ be  an $R$-module. If $M$ is $\omega$-$\infty$-torsionfree as an $R$-module, then $S\otimes_RM$ is $({S\otimes_R\omega})$-$\infty$-torsionfree as an $S$-module.
\end{cor}

\begin{cor}{\label{torsionfree equ II}} Let $S/R$ be a Frobenius extension  with $_RS_R$ centrally projective over $R$ and $M$ be an $S$-module. Then  $M$ is $({_SS}\otimes_R\omega)$-$\infty$-torsionfree as an $S$-module if and only if $M$ is $\omega$-$\infty$-torsionfree as an $R$-module.
\end{cor}

\begin{lem}{\label{orthonogal equivalence I}}Let $S/R$  be a Frobenius extension and $M$ be an $S$-module. Then $M\in{^\bot_R \omega}$ as an $R$-module if and only if $M\in{^\bot(_SS\otimes_R\omega)}$ as an $S$-module.
\end{lem}
\noindent{\bf Proof.} Since $S/R$ is a Frobenius extension, we have $S\otimes_RM\cong {\rm Hom}_R({_RS_S,{_RM}})$ for any
${_RM}\in R$-${\rm mod}$.
By the adjoint isomorphism, for any $i\geq 0$, we have
\begin{align*}{\rm Ext}^i_R(M,\omega)
&\cong{\rm Ext}^{i}_R({_RS\otimes_SM},\omega)\\
&\cong{\rm Ext}^{i}_S({_SM}, {\rm Hom}_R({_RS_S}, \omega))\\
&\cong {\rm Ext}^i_S({_SM}, {_SS\otimes_R\omega}).
\end{align*}
Consequently, ${\rm Ext}^i_R({_RM},\omega)=0$ for any $i\geq 1$ if and only if ${\rm Ext}^i_S({_SM}, {_SS\otimes_R\omega})=0$ for any $i\geq 1$. \hfill$\square$

\begin{lem}{\label{orthonogal equivalence II}}Let $S/R$  be a Frobenius extension with $_RS_R$ centrally projective over $R$  and $M$ be an $R$-module. If $M\in{^\bot \omega}$, then $_SS\otimes_RM\in{^\bot(_SS\otimes_R\omega)}$.
\end{lem}

\noindent{\bf Proof.} By the adjoint isomorphism, for any $i\geq 1$, we have
\begin{align*}
{\rm Ext}^{i}_S(S\otimes_RM,_SS\otimes_R\omega)
&\cong{\rm Ext}^{i}_R(M, {\rm Hom}_S({_SS_R}, _SS\otimes_R\omega))\\
&\cong {\rm Ext}^i_R(M, _RS\otimes_R\omega).
\end{align*}

Since $S/R$ is a Frobenius extension and $_RS_R$ is centrally projective over $R$, we know that $_RS\otimes_R\omega\in{\rm add}{_R\omega}$.
 By assumption, ${_RM}\in{^\bot\omega}$, and so ${\rm Ext}^i_R(M,  {_RS\otimes_R\omega})=0$  for any $i\geq 1$.
It follows that ${\rm Ext}^{i}_S(S\otimes_RM, {_SS\otimes_R\omega})=0$ for any $i\geq 1$, that is, $S\otimes_RM\in {^\bot{_S(S\otimes_R\omega)}}$.
\hfill$\square$

Combining Corollary \ref{torsionfree equ II} and Lemma \ref{orthonogal equivalence I}, we have
\begin{thm} {\label{equivalent of G-dimension}}Let $S/R$  be a Frobenius extension with $_RS_R$ centrally projective over $R$  and $M$ be an $S$-module. Then $M$ is of generalized G-dimension zero with respect to $_SS\otimes_R\omega$ as an $S$-module if and only if $M$ is  of generalized G-dimension zero with respect to $_R\omega$ as an $R$-module.
\end{thm}
\noindent{\bf Proof.} By Lemma \ref{lem G G-dimension zero}, $M$ is of generalized G-dimension zero with respect to $_SS\otimes_R\omega$ as an $S$-module if and only if $M\in{^\bot_S(S\otimes_R\omega)}$ and $M$ is $(S\otimes_R\omega)$-$\infty$-torsionfree, which is equivalent to  $M\in{^\bot_R\omega}$ and $M$ is $\omega$-$\infty$-torsionfree as an $R$-module by  Lemma \ref{orthonogal equivalence I} and Corollary \ref{torsionfree equ II}, respectively. The latter condition is equivalent to  $M$ being of generalized G-dimension zero with respect to $_R\omega$ as an $R$-module. \hfill$\square$

Similarly, by Corollary \ref{torsionfree equ I} and Lemma \ref{orthonogal equivalence II}, we have
\begin{prop}{\label{equivalent of G-dimension zero II}} Let $S/R$  be a Frobenius extension with $_RS_R$ centrally projective over $R$  and $M$ be an $R$-module. If $M$ is an $R$-module of generalized G-dimension zero with respect to $\omega$, then  ${S\otimes_RM}$ is an $S$-module of generalized G-dimension zero with respect to ${S\otimes_R\omega}$.
\end{prop}

 For a module $M\in R$-mod, the \emph{generalized G-dimension} of
$M$ with respect to $\omega$, denoted by G-dim$_\omega(M)$, is defined as G-dim$_\omega(M)={\rm inf}\{n\mid \exists$ exact sequence $0\rightarrow G_n\rightarrow \cdots\rightarrow G_1\rightarrow G_0\rightarrow M\rightarrow 0 $ with $G_i\in \mathcal{G}_\omega$ for $0\leq i\leq n\}$.   We have
G-dim$_\omega(M)\geq 0$, and we set G-dim$_\omega(M)=\infty$ if no such integer exists.

Recall that a subcategory $\mathcal{X}$ of $R$-mod is said to be \emph{projectively resolving} if it contains all projective modules, and for every short exact sequence $0\rightarrow X'\rightarrow X\rightarrow X''\rightarrow 0$ with $X''\in\mathcal{X}$, $X'\in\mathcal{X}$ if and only if $X\in\mathcal{X}$. By  \cite[Lemmas 5.7-5.8]{Hu3}, the subcategory $\mathcal{G}_\omega$  is projectively resolving. The following is a consequence of Theorem 3.12 in \cite{AB}.

\begin{lem}\label{description of G-dim} For an integer $n\geq 0$ and a module $M$ in $R$-$\rm mod$, the following statements are equivalent.

$(1)$ $\rm G$-${\rm dim}_\omega(M)\leq n$;

$(2)$ If there exists an exact sequence $0\rightarrow G_n\rightarrow G_{n-1}\rightarrow \cdots \rightarrow G_1\rightarrow G_0\rightarrow M\rightarrow 0$  with $G_i$ of generalized $\rm G$-dimension zero with respect to $\omega$ for $0\leq i\leq n-1$, then $G_n$ is also a module of generalized $\rm G$-dimension zero with respect to $\omega$.
\end{lem}

The following shows that the generalized G-dimension with respect to a Wakamatsu tilting module  is preserved under Frobenius extensions.

\begin{thm}\label{equivalent of generlazied G-dim} Let $S/R$ be a Frobenius extension with $_RS_R$ centrally projective over $R$ and $M$ be an $S$-module. Then $\rm G$-${\rm dim}_{S\otimes_R\omega}(M)$=$\rm G$-${\rm dim}_\omega(M)$.
\end{thm}
\noindent{\bf Proof.} Without loss of generality, we assume that G-${\rm dim}_{S\otimes_R\omega}(M)=n<\infty$. There is an exact sequence
$0\rightarrow G_n\rightarrow \cdots \rightarrow G_1\rightarrow G_0\rightarrow {_SM}\rightarrow 0 $ in $S$-mod such that $G_i$ is of generalized G-dimension zero with respect to ${S\otimes_R\omega}$ for $0\leq i\leq n$.
Applying  the restriction functor, we get the  exact sequence
$0\rightarrow G_n\rightarrow \cdots\rightarrow G_1\rightarrow G_0\rightarrow {_RM}\rightarrow 0 $ in $R$-mod with  $G_i$ being of generalized G-dimension zero with respect to $\omega$ for $0\leq i\leq n$ by Theorem \ref{equivalent of G-dimension}. Hence G-${\rm dim}_\omega(M)\leq n=G$-${\rm dim}_{{S\otimes_R\omega}}(M)$.

Conversely, we can assume that G-${\rm dim}_\omega(M)=m<\infty$. Regarding $M$ as an $S$-module, we choose  an exact sequence
$0\rightarrow K_m\rightarrow G_{m-1}\rightarrow\cdots\rightarrow G_1\rightarrow G_0\rightarrow {_SM}\rightarrow 0$ in $S$-mod such that $G_i$ is of generalized  G-dimension zero with respect to ${S\otimes_R\omega}$ for $0\leq i\leq {m-1}$ (in fact, the  $G_i$ may all be chosen to be  projective $S$-modules).
Applying the restriction functor, we get the exact sequence $0\rightarrow K_m\rightarrow G_{m-1}\rightarrow\cdots\rightarrow G_1\rightarrow G_0\rightarrow {_RM}\rightarrow 0$ in $R$-mod  such that $G_i$ is of generalized G-dimension zero with respect to $\omega$  for $0\leq i\leq {m-1}$ by Theorem \ref{equivalent of G-dimension}.
Since G-${\rm dim}_\omega(M)=m$, $K_m$ is also of generalized G-dimension zero with respect to $\omega$ as an $R$-module by Lemma \ref{description of G-dim}.
Again by Theorem \ref{equivalent of G-dimension}, $K_m$ is also of generalized G-dimension zero with respect to ${S\otimes_R\omega}$ as an $S$-module. Then G-${\rm dim}_{{S\otimes_R\omega}}(M)\leq m=$G-${\rm dim}_\omega(M)$.
\hfill$\square$

Using Proposition \ref{equivalent of G-dimension zero II} and an argument similar to that of Theorem \ref{equivalent of generlazied G-dim}, we obtain the following result.
\begin{prop}\label{equivalent of G-dim II} Let $S/R$ be a Frobenius extension with $_RS_R$ centrally projective over $R$ and $M$ be an $R$-module. Then $\rm G$-${\rm dim}_{S\otimes_R\omega}(S\otimes_RM)\leq{\rm G}$-${\rm dim}_\omega(M)$.
\end{prop}
\noindent{\bf Proof.} Without loss of generality, we assume that G-${\rm dim}_{\omega}(M)=n<\infty$. There is an exact sequence
$0\rightarrow G_n\rightarrow \cdots\rightarrow G_1\rightarrow G_0\rightarrow {M}\rightarrow 0 $ in $R$-mod such that $G_i$ is of generalized G-dimension zero with respect to $\omega$ for $0\leq i\leq n$.
Applying  the exact functor $_SS\otimes_R-$, we get the  exact sequence
$0\rightarrow S\otimes_RG_n\rightarrow \cdots\rightarrow S\otimes_RG_1\rightarrow S\otimes_RG_0\rightarrow {_SS}\otimes_RM\rightarrow 0 $ in $S$-mod with  $S\otimes_RG_i$ having generalized G-dimension zero with respect to $S\otimes_R\omega$ for $0\leq i\leq n$ by Proposition \ref{equivalent of G-dimension zero II}. Hence G-${\rm dim}_{S\otimes_R\omega}(S\otimes_RM)\leq n=G$-${\rm dim}_\omega(M)$.                      \hfill$\square$

When $\omega=R$, we have
${}_R(S\otimes_RR)\cong{}_RS\in\operatorname{add}({}_RR)$,
since $S/R$ is a Frobenius extension.
Thus the central projectivity assumption is unnecessary
in this case, and the proof of Theorem~5.9 yields
the following corollary, which is one of the main results in \cite{R,Z}.

\begin{cor} Let $S/R$ be a Frobenius extension and $M$ an $S$-module. Then we have {\rm Gpd}$_S(M)$={\rm Gpd}$_R(M)$, where {\rm Gpd}$_S(M)$ denotes the Gorenstein projective dimension of $_SM$.
\end{cor}

\noindent{\bf Acknowledgements}

The authors thank Professor Xiao-Wu Chen for his helpful suggestions. The authors are very grateful to the referee for many helpful comments and suggestions. This work is supported by the National Natural Science Foundation of China (No. 12371015 and 12571034) and the Excellent University Research and Innovation Team in Anhui Province (No. 2024AH010002).

\vspace{0.5cm}

 {\footnotesize \noindent  Yanhong Bao, Zhibing Zhao\\
 Center for Pure Mathematics, School of Mathematical Sciences, Anhui University,  Hefei 230601, Anhui, PR China\\

\footnotesize \noindent Jiafeng L$\ddot{u}$\\
School of Mathematical Sciences, Zhejiang Normal University, Jinhua 321004, Zhejiang, PR China}

\end{document}